\documentclass{commat}

\usepackage{dsfont}
\usepackage{enumerate}
\usepackage{longtable}

\title{%
   Cyclicity of the $2$-class group of the first Hilbert $2$-class field of some number fields
    }

\author{%
    A. Azizi, M. Rezzougui and A. Zekhnini
    }

\affiliation{
    \address{Abdelmalek Azizi --
    Mohammed First University, Sciences Faculty, Mathematics Department, Oujda, Morocco
        }
    \email{
    abdelmalekazizi@yahoo.fr
    }
    \address{Mohammed Rezzougui --
    Mohammed First University, Sciences Faculty, Mathematics Department, Oujda, Morocco
        }
    \email{
    morez2100@hotmail.fr
    }
    \address{Abdelkader Zekhnini --
    Mohammed First University, Sciences Faculty, Mathematics Department, Oujda, Morocco
        }
    \email{
    zekha1@yahoo.fr
    }
    }

\abstract{%
     Let $\mathds{k}$ be a real quadratic number field. Denote by $\mathrm{Cl}_2(\mathds{k})$ its $2$-class group and by $\mathds{k}_2^{(1)}$ (resp. $\mathds{k}_2^{(2)}$) its first (resp. second) Hilbert $2$-class field. The aim of this paper is
to  study, for a real quadratic number field whose discriminant is divisible by one prime number congruent to $3$ modulo 4, the metacyclicity  of $G=\mathrm{Gal}(\mathds{k}_2^{(2)}/\mathds{k})$ and the cyclicity of   $\mathrm{Gal}(\mathds{k}_2^{(2)}/\mathds{k}_2^{(1)})$ whenever the rank of $\mathrm{Cl}_2(\mathds{k})$ is $2$, and the $4$-rank of $\mathrm{Cl}_2(\mathds{k})$ is $1$.
    }

\keywords{%
   quadratic field, Hilbert $2$-class field, 2-class group, metacyclic 2-group, metabelian 2-group
    }

\msc{%
   Primary 11R29,  11R11, 11R20, 11R32, 11R37; Secondary 20D15
    }

\VOLUME{32}
\NUMBER{1}
\YEAR{2024}
\firstpage{157}
\DOI{https://doi.org/10.46298/cm.10983}

\begin{paper}

\section{Introduction}
Let $\mathds{k}$ be  an algebraic number field and $\mathrm{Cl}_2(\mathds{k})$ its $2$-class group, that is, the $2$-Sylow subgroup of its ideal class group $\mathrm{Cl}(\mathds{k})$. Let  $\mathds{k}^{(1)}_2$ be the Hilbert $2$-class field of $\mathds{k}$, that is, the maximal abelian extension of $\mathds{k}$ everywhere
unramified of $2$-power degree over $\mathds{k}$. Put $\mathds{k}^{(0)}_2 = \mathds{k}$ and let $\mathds{k}^{(i+1)}_2$ denote the Hilbert $2$-class field of $\mathds{k}_2^{(i)}$
for any integer $i\geq 0$. Then the sequence of fields
$$\mathds{k}=\mathds{k}^{(0)}_2 \subseteq  \mathds{k}^{(1)}_2 \subseteq   \mathds{k}^{(2)}_2  \subseteq  \cdots \subseteq   \mathds{k}^{(i)}_2\subseteq  \cdots, $$
is called   the $2$-class field tower of $\mathds{k}$.

If, for all $i\geq0$, we have $\mathds{k}_2^{(i)}\neq \mathds{k}_2^{(i+1)}$, the tower is said to be infinite. Otherwise, the tower is said to be  finite, and the minimal integer $i$ such that $\mathds{k}_2^{(i)}= \mathds{k}_2^{(i+1)}$ is called the length of the tower.
Unfortunately, there is no known method to decide  whether or not a $2$-class field tower of a number field  is infinite. However,  it is known that if the rank of $\mathrm{Cl}_2(\mathds{k}^{(1)}_2)$ is at most $2$, then the tower is finite and its length  is at most $3$ (cf. \cite{Bl-57}); it is also known that if the rank of $\mathrm{Cl}_2(\mathds{k}^{(1)}_2)$ is equal to $3$, then there are fields  $\mathds{k}$ with infinite $2$-class field tower (cf. \cite{Sch80}). Therefore, it is  interesting  to determine all fields such that rank$(\mathrm{Cl}_2(\mathds{k}^{(1)}_2))\leq2$. That is why Benjamin et al. started a project which aims to characterize all quadratic fields $\mathds{k}$  satisfying the last  condition  (cf. \cite{Be-17}, \cite{B.L.S-97}, \cite{B.L.S-98}, \cite{BLS-01}, \cite{BS-17}, \cite{BS-18}, \cite{BS-182}, \cite{CD-92}).  Our present paper as well as our previous one (see \cite{AMRZ19}) are  part of this project.

We aim to  study the cyclicity of $\mathrm{Cl}_2(\mathds{k}^{(1)}_2)$ of  real quadratic fields $\mathds{k}$ such that $\mathrm{Cl}_2(\mathds{k})$ is of the form $(2^n, 2^m)$ for some $n\geq1$ and $m\geq2$, and their discriminants $d_\mathds{k}$ are divisible by  primes congruent to $3$ modulo $4$. In this paper, which is a continuation of \cite{AMRZ19}, we consider the  field $\mathds{\mathds{k}}=\mathbb{Q}(\sqrt{2p_1p_2q})$, where $p_1\equiv p_2\equiv-q\equiv 1 \pmod4$ are primes and $\mathrm{Cl}_2(\mathds{k})\simeq (2, 2^n)$, with $n\geq2$. We determine complete criteria for $G=\mathrm{Gal}(\mathds{k}_2^{(2)}/\mathds{k})$ to be metacyclic and complete criteria for $\mathrm{Cl}_2(\mathds{k}^{(1)}_2)$ to be cyclic whenever $G$ is not metacyclic.

\section{Preliminary results}

We begin by collecting some results that will be useful later. We recall that a $2$-group $G$ is said to be metacyclic if there exists a normal cyclic subgroup $N$ of $G$ such that  $G/N$ is cyclic. It is known that if  $G$ is metacyclic, then the minimal  number of its generators is less or equal to $2$; this number is called the rank of $G$ and will be denoted by $d(G)$. On the other hand, if $d(G)=2$, then $G/G'$ is of  type $(2^n, 2^m)$ with  $n$ and $m\in \mathds{N}^*$, where $G'$ is the commutator subgroup of $G$. If $n=m=1$, then it is known  that $G$ is dihedral,  semi-dihedral, quaternionic or abelian of type $(2, 2)$ (cf.\ \cite{Ki76}, \cite{Go68}). In these cases, $G$ admits a cyclic maximal subgroup, and thus is metacyclic. 
 By Blackburn \cite{Bl-61}, we know  that the metacyclicity of a $2$-group $G$ is characterized by the rank of its maximal subgroups, and we have the following lemmas.

\begin{lemma}[\cite{AMRZ19}]\label{lem7}
 Let $G$ be a finite  $2$-group   such that  $G/G'$ is of type $(2^n, 2^m)$, where $n\geq1$ and $m\geq2$. Denote by $H_i$ $(i=1, 2, 3)$, the three maximal subgroups of  $G$. Then
$G$  is metacyclic    if and only if  $d(H_i)\leq 2$ for all $ i=1, 2, 3$.
\end{lemma}
 \begin{lemma}[\cite{B.L.S-97}] \label{lem6}  Let $G$ be a non-metacyclic  $2$-group   such that   $G/G'$ is of type $(2, 2^m)$, where $m\geq2$. Then $G$ admits two maximal subgroups  $H_1$ and $H_2$ such that $H_1/G'$ and $H_2/G'$ are cyclic. Moreover, if $G'$ is cyclic, then $H_1$ and $H_2$ are metacyclic.
\end{lemma}

We continue by fixing some notation. For a number field $\mathds{k}$, denote by $\mathrm{Cl}_2(\mathds{k})$
 its $2$-class group in the ordinary sense, denote by $h_2(\mathds{k})$ the order of $\mathrm{Cl}_2(\mathds{k})$, denote by $\mathds{k}^{(1)}_2$  the
Hilbert $2$-class field of $\mathds{k}$, and denote by $\mathds{k}^{(2)}_2$ the Hilbert $2$-class field of $\mathds{k}^{(1)}_2$. If $G = \mathrm{Gal}(\mathds{k}^{(2)}_2/\mathds{k})$, then it
is well known from class field theory that $G' = \mathrm{Gal}(\mathds{k}^{(2)}_2/\mathds{k}^{(1)}_2) \simeq \mathrm{Cl}_2(\mathds{k}^{(1)}_2)$ and $G/G' = \mathrm{Gal}(\mathds{k}^{(1)}_2/\mathds{k}) \simeq \mathrm{Cl}_2(\mathds{k})$.
 Note that if $\mathrm{Cl}_2(\mathds{k})$ is of type $(2, 2^n)$, with $n\geq2$, then $\mathds{k}$ admits three unramified quadratic extensions within $\mathds{k}^{(1)}_2$, which will be denoted by $\mathbb{K}_i$ ($i= 1, 2, 3$). We suppose that $\mathbb{K}_3$ is included in the three unramified biquadratic
extensions of $\mathds{k}$ within $\mathds{k}^{(1)}_2$. The following result was shown in our earlier paper \cite{ARZ20}.
\begin{theorem}\label{thm0}
Keep the notations above and assume  $G/G'$ is of type $(2, 2^n)$, where $n\geq2$. Then
\begin{enumerate}[\rm1.]
\item $G$ is abelian or modular  if and only if
$$\mathrm{rank}(\mathrm{Cl}_2(\mathbb{K}_i))=1\;(i= 1, 2)\;\; \text{and}\;\; \mathrm{rank}(\mathrm{Cl}_2(\mathbb{K}_3))=2.$$
 \item  $G$  is metacyclic non-abelian non-modular if and only if  $$\mathrm{rank}(\mathrm{Cl}_2(\mathbb{K}_i))=2\;\; \text{for all}\;\; i =1, 2, 3.$$
 \item $G$  is non-metacyclic non-abelian  if and only if
$$\mathrm{rank}(\mathrm{Cl}_2(\mathbb{K}_i))=2\;(i= 1, 2)\;\; \text{and}\;\; \mathrm{rank}(\mathrm{Cl}_2(\mathbb{K}_3))=3.$$
 \end{enumerate}
\end{theorem}

Let $\mathds{k} = \mathbb{Q}(\sqrt d)$ be an arbitrary quadratic number field with a square-free integer $d$, and $d_\mathds{k}$ be its  discriminant. For a prime number $p$, define:
 \[
 p^*=
 \begin{cases}
 (-1)^{\frac{p-1}{2}} p, & \textup{if } p \ne 2; \\
 -4, & \textup{if } p = 2 \textup{ and } d \equiv 3 \textup{ (mod $4$)}; \\
 8, & \textup{if } p = 2 \textup{ and } d \equiv 2 \textup{ (mod $8$)}; \\
 -8, & \textup{if } p = 2 \textup{ and } d \equiv -2 \textup{ (mod $8$)}.
 \end{cases}
 \]
 Then, let $d_\mathds{k}=p^*_1 \ldots p^*_sp^*_{s+1}\ldots p^*_{s+t} $  such that $p^*_1, \ldots, p^*_s$ are positive and $p^*_{s+1},\ldots, p^*_{s+t}$ are negative. The Rédei matrix $R_\mathds{k}$ is defined to be the matrix in $M_{(s+t)\times (s+t)}(\mathbb{Z}/2\mathbb{Z})$ with entries $a_{i,j}$ given by:
$(-1)^{a_{i,j}}= \left( \frac{p^*_i}{p_j} \right)$ if $i\neq j$ and
$(-1)^{a_{i,j}}=\left( \frac{d_\mathds{k}/p^*_i}{p_i} \right)$ if  $i= j$,
where $\left( \frac{\bullet}{\bullet} \right)$ is the  Legendre symbol. Then the $4$-rank of $\mathrm{Cl}^{+}(\mathds{k})$,  the class group of $\mathds{k}$ in  the  narrow sense, is given by:
\begin{theorem}[\cite{LR}]\label{1} Let $\mathds{k}$ be a quadratic number field, then
$$ 4\text{-}\mathrm{rank}(\mathrm{Cl}^{+}(\mathds{k}))=s+t-1- \mathrm{rank}(R_\mathds{k}).$$
\end{theorem}
\begin{remark}\label{rema1} If $d_\mathds{k}$ is divisible by  a prime congruent to   $3$ modulo $4$, then
\[
\mathrm{Cl}_2^{+}(\mathds{k}))\simeq \mathbb{Z}/2\mathbb{Z}\times \mathrm{Cl}_2(\mathds{k})
\qquad \textup{and} \qquad
4\text{-}\mathrm{rank}(\mathrm{Cl}^{+}(\mathds{k}))=4\text{-}\mathrm{rank}(\mathrm{Cl}(\mathds{k})),
\]
where $\mathrm{Cl}_2^{+}(\mathds{k})$ is  the $2$-class group of $k$ in  the  narrow sense.
\end{remark}

We make use of the well known Kuroda Class Number Formula, which we state as the following theorem.

\begin{theorem} [\cite{Lm}] \label{thm5}Let $ \mathbb{K}/\mathds{k}$ be an arbitrary normal quartic extension of number fields with Galois group of type $(2,2)$, and let $\mathbb{K}_j$ $(j=1, 2, 3)$ denote the quadratic subextensions. Then the class number of $\mathbb{K}$ satisfies
$$ h(\mathbb{K})= 2^{d-\kappa-2-\upsilon} \frac{q(\mathbb{K}/\mathds{k}) h(\mathbb{K}_1)h(\mathbb{K}_2)h(\mathbb{K}_3)}{(h(\mathds{k}))^2},$$
where $q(\mathbb{K}/\mathds{k})=[E_\mathbb{K}:E_1E_2 E_3]$ denotes the unit index of  $ \mathbb{K}/\mathds{k}\ ($with $E_j=$ the unit group of $\mathbb{K}_j)$, $d$ is the number of infinite primes in $\mathds{k}$ that ramify in $\mathbb{K}$, $\kappa$ is the $\mathbb{Z}$-rank of the unit group $E_\mathds{k}$ of $\mathds{k}$, and $\upsilon=0$ except when $\mathbb{K} \subseteq \mathds{k}(\sqrt{E_\mathds{k}})$, in which case  $\upsilon=1$.
\end{theorem}

To prove our main theorems, we also need the following results.
\begin{theorem}[\cite{Be-17}]\label{thm4} Let $\mathds{k}$  be a number field such that  $\mathrm{Cl}_2(\mathds{k})\simeq (2, 2^n)$, where $n\geq 2$. Denote by  $\mathbb{K}_i$ $(i=1, 2, 3)$,  the three unramified quadratic extensions of $\mathds{k}$. Then the  $2$-class group of  $\mathds{k}_2^{(1)}$ is a non-elementary  cyclic group if and only if  $h_2(\mathbb{K}_i)\geq 2h_2(\mathds{k})$ and $h_2(\mathbb{K}_j)= h_2(\mathbb{K}_m)=h_2(\mathds{k})$ for some $\{i, j, m \}=\{1, 2, 3\}$.
\end{theorem}

\begin{lemma}[\cite{Be-17}] \label{Lem5}Let $\mathds{k}$  be a number field such that  $\mathrm{Cl}_2(\mathds{k})\simeq (2^m, 2^n)$, $m\geq1$, $n\geq1$. Denote by   $\mathbb{K}_i$ $(i=1, 2, 3)$,  the three unramified quadratic extensions of $\mathds{k}$. Then
$h_2(\mathds{k}_2^{(1)})=\nolinebreak2$ if and only if $h_2(\mathbb{K}_0 )= (1/2) h_2(\mathds{k})$ where $\mathbb{K}_0=\mathbb{K}_1\mathbb{K}_2\mathbb{K}_3$.
\end{lemma}

\begin{corollary}[\cite{Be-17}] \label{coro5}Let $\mathds{k}$  be a real quadratic number field such that  $ \mathrm{Cl}_2(\mathds{k})\simeq (2^m, 2^n)$, $m\geq1$, $n\geq1$. Denote by  $\mathbb{K}_i$ $(i=1, 2, 3)$,  the three unramified quadratic extensions of $\mathds{k}$. Then  $h_2(\mathds{k}_2^{(1)})=2$ if and only if  $h_2(\mathbb{K}_1)=h_2(\mathbb{K}_2)= h_2(\mathbb{K}_3)=h_2(\mathds{k})$ and $q(\mathbb{K}_0/\mathds{k})=4$, where $\mathbb{K}_0=\mathbb{K}_1\mathbb{K}_2\mathbb{K}_3$.
\end{corollary}
\begin{theorem}[\cite{Be-17}]\label{thm7} Let $\mathds{k}$ be a real quadratic number field with  $\mathrm{Cl}_2(\mathds{k})\simeq (2^m, 2^n)$, $m\geq1$, $n\geq2$, and $d_\mathds{k}=d_1d_2r_1r_2$ or $r_1r_2r_3r_4$ be its discriminant, where $d_1$ and $d_2$ are positive prime discriminants and $r_1$, $r_2$, $r_3$, $r_4$ are negative prime discriminants. Denote by  $\mathbb{K}_i$ $(i=\nolinebreak1, 2, 3)$,  the three unramified quadratic extensions of $k$. If $h_2(\mathds{k}_2^{(1)})=2$ then $Q_{\mathbb{K}_i}=\nolinebreak Q_{\mathbb{K}_j}=2$, and $Q_{\mathbb{K}_s}=1$ or $2$ for some $\{i, j, s\}=\{1, 2, 3\}$, where $Q_{\mathbb{K}}$ denotes the unit index of $\mathbb{K}$.
\end{theorem}

\section{The $4$-rank of the $2$-class group of $\mathds{k}=\mathbb{Q}(\sqrt{2p_1p_2q})$.}

Let $p_1\equiv p_2\equiv-q\equiv 1 \!\!\pmod4$ be different positive prime integers and $\mathds{k}=\mathbb{Q}(\sqrt{2p_1p_2q})$. It is well known, by genus theory, that the $2$-rank of the class group of $\mathds{k}$ is $2$. The purpose of this section is to determine the $4$-rank of the $2$-class group of $\mathds{k}$.

 E. Benjamin and C. Snyder characterized  real quadratic fields whose $2$-class group is  of type   $(2, 2)$ in \cite{Be-Sn}. In particular, they proved the following theorem.
\begin{theorem}\label{thm2}  Let $p_1\equiv p_2\equiv-q\equiv 1 \pmod4$ be different prime integers. Then the $2$-class group of $ \mathds{k}=\mathbb{Q}(\sqrt{2p_1p_2q} )$  is of type $(2, 2)$ (i.e.\ $4\text{-}\mathrm{rank}(\mathrm{Cl}_2(\mathds{k}))=0)$ if and only if  one of the following conditions is satisfied.
\begin{enumerate}[\rm1.]
\item 
\begin{itemize}
    \item $\left( \frac{p_1}{p_2} \right)=1$, and 
    \item either $\left( \frac{2}{p_1} \right)=-1$ or $\left( \frac{q}{p_1} \right)=-1$, and 
    \item either $\left( \frac{2}{p_2} \right)=-1$ or $\left( \frac{q}{p_2} \right)=-1$, and
    \item $\left( \frac{2}{p_1} \right)$, $\left( \frac{2}{p_2} \right)$, $\left( \frac{q}{p_1} \right)$, $\left( \frac{q}{p_2} \right)$ are not all equal.
\end{itemize}
\item $\left( \frac{p_1}{p_2} \right)=-1$ and  $\left( \frac{2}{p_1} \right)$, $\left( \frac{2}{p_2} \right)$, $\left( \frac{q}{p_1} \right)$, $\left( \frac{q}{p_2} \right)$ are not all equal.
\end{enumerate}
\end{theorem}

In the following theorem, we give necessary and sufficient conditions for the  $2$-class group of $ \mathds{k}=\mathbb{Q}(\sqrt{2p_1p_2q} )$  to be of type $(2, 2^n)$ or $(2^m, 2^n)$, where $n\geq2$ and $m\geq2$.
\begin{theorem} \label{thm3}
 Let $p_1\equiv p_2\equiv-q\equiv 1 \pmod4$ be different prime integers. Then the $2$-class group of $ \mathds{k}=\mathbb{Q}(\sqrt{2p_1p_2q} )$  is of type $(2,2^n)$, where $n\geq2,  $(i.e.\  $ \text{$4$-rank}(\mathrm{Cl}_2(\mathds{k}))=1$) if and only if  one of the following conditions is satisfied.
\begin{enumerate}[\rm1.]
\item  $\left( \frac{2}{p_1} \right)=\left( \frac{2}{p_2} \right)=\left( \frac{p_1}{p_2} \right)=1$ and $\left( \frac{q}{p_1} \right) \left( \frac{q}{p_2} \right)=-1$.
\item  $\left( \frac{2}{p_1} \right)=\left( \frac{2}{p_2} \right)=\left( \frac{q}{p_1} \right)=\left( \frac{q}{p_2} \right)=1$ and $\left( \frac{p_1}{p_2} \right)=-1$.
\item  $\left( \frac{2}{p_1} \right)=-\left( \frac{2}{p_2} \right)=1$ and $\left( \frac{p_1}{p_2} \right)=\left( \frac{q}{p_1} \right)=1$.
\item $\left( \frac{2}{p_1} \right)=\left( \frac{2}{p_2} \right)=\left( \frac{q}{p_1} \right)=\left( \frac{q}{p_2} \right)=-1$.
\end{enumerate}
Moreover,   $4\text{-}\mathrm{rank}(\mathrm{Cl}_2(\mathds{k}))=2$ if and only if $$\left( \frac{2}{p_1} \right)=\left( \frac{2}{p_2} \right)=\left( \frac{p_1}{p_2} \right)=\left( \frac{q}{p_1} \right)=\left( \frac{q}{p_2} \right)=1.$$
\end{theorem}
\begin{proof} Proceeding as in \cite{AMRZ19}, the results are deduced by applying Theorem \ref{1} and Remark~\ref{rema1}.
\end{proof}
\section{The $\mathrm{FSU}$s of certain biquadratic number fields}
  Let $p_1\equiv p_2\equiv-q\equiv 1 \pmod4$ be different prime integers. Put $\mathds{k}=\mathbb{Q}(\sqrt{2p_1p_2q})$. Consider the following three unramified quadratic extensions of $\mathds{k}$:
  \begin{center}$\mathbb{K}_1=\mathbb{Q}( \sqrt{p_1},\sqrt{2p_2q} )$,\quad $ \mathbb{K}_2=\mathbb{Q}( \sqrt{p_2},\sqrt{2p_1q} )$\  and\   $ \mathbb{K}_3=\mathbb{Q}( \sqrt{2q},\sqrt{p_1p_2} )$.\end{center}
   Let  $\varepsilon_{2p_1p_2q}=x+y\sqrt{2p_1p_2q}$, $\varepsilon_{2p_1q}=z+t\sqrt{2p_1q}$ and $\varepsilon_{2p_2q}=a+b\sqrt{2p_2q}$ be the fundamental units of  $\mathbb{Q}(\sqrt{2p_1p_2q})$, $\mathbb{Q}(\sqrt{2p_2q})$ and $\mathbb{Q}(\sqrt{2p_1q})$ respectively. The goal of this section is to determine a Fundamental System of Units  ($\mathrm{FSU}$) of $\mathbb{K}_i$  basing on  the conditions cited in Theorem \ref{thm3}.

    Using similar arguments as in the proof of Lemma 4.1 of \cite{AMRZ19} (see also \cite{ATZAR}), we get the following lemmas.
\begin{lemma}\label{lem2} Suppose that  $\left( \frac{2}{p_1} \right)=\left( \frac{2}{p_2} \right)=\left( \frac{p_1}{p_2} \right)=\left( \frac{q}{p_1} \right)=-\left( \frac{q}{p_2} \right)=1$.
\begin{enumerate}[\rm1.]
\item If $x\pm1$   is a square in $\mathbb{N}$,  then
\begin{enumerate}[\rm i.]
\item  $\{\varepsilon_{p_1}, \varepsilon_{2p_2q}, \sqrt{\varepsilon_{2p_2q}\varepsilon_{2p_1p_2q}} \}$  is a  $\mathrm{FSU}$ of  $\mathbb{K}_1$.
\item  $\{\varepsilon_{p_2}, \varepsilon_{2p_1q}, \sqrt{\varepsilon_{2p_1q}\varepsilon_{2p_1p_2q}} \}$ or $\{\varepsilon_{p_2}, \varepsilon_{2p_1q}, \varepsilon_{2p_1p_2q} \}$  is a  $\mathrm{FSU}$ of  $\mathbb{K}_2$ according as  $z\pm1$   is or not a square in $\mathbb{N}$.
\item $\{\varepsilon_{2q}, \varepsilon_{p_1p_2}, \sqrt{\varepsilon_{2q}\varepsilon_{2p_1p_2q}} \}$  is a  $\mathrm{FSU}$ of  $\mathbb{K}_3$.
\end{enumerate}
\item If $p_1(x\pm1)$   is a square in $\mathbb{N}$,  then
\begin{enumerate}[\rm i.]
\item $\{\varepsilon_{p_1}, \varepsilon_{2p_2q}, \sqrt{\varepsilon_{2p_2q}\varepsilon_{2p_1p_2q}} \}$  is a  $\mathrm{FSU}$ of  $\mathbb{K}_1$.
\item $\{\varepsilon_{p_2}, \varepsilon_{2p_1q}, \sqrt{\varepsilon_{2p_1q}\varepsilon_{2p_1p_2q}} \}$ or $\{\varepsilon_{p_2}, \varepsilon_{2p_1q}, \varepsilon_{2p_1p_2q} \}$  is a  $\mathrm{FSU}$ of  $\mathbb{K}_2$ according as  $p_1(z\pm1)$   is or not a square in $\mathbb{N}$.
\item $\{\varepsilon_{2q},\, \varepsilon_{p_1p_2},\, \sqrt{\varepsilon_{2q}\varepsilon_{p_1p_2}\varepsilon_{2p_1p_2q}} \}$ or  $\{\varepsilon_{2q},\, \varepsilon_{p_1p_2},\, \varepsilon_{2p_1p_2q} \}$  is a  $\mathrm{FSU}$ of  $\mathbb{K}_3$ according to whether $N_{\mathbb{Q}( \sqrt{p_1p_2} )/\mathbb{Q}}(\varepsilon_{p_1p_2})$ equals $1$ or $-1$.
\end{enumerate}
\item If $2p_1(x\pm1)$   is a square in $\mathbb{N}$,  then
\begin{enumerate}[\rm i.]
\item $\{\varepsilon_{p_1}, \varepsilon_{2p_2q}, \sqrt{\varepsilon_{2p_1p_2q}} \}$  is a  $\mathrm{FSU}$ of  $\mathbb{K}_1$.
\item $\{\varepsilon_{p_2}, \varepsilon_{2p_1q}, \sqrt{\varepsilon_{2p_1q}\varepsilon_{2p_1p_2q}} \}$ or $\{\varepsilon_{p_2}, \varepsilon_{2p_1q}, \varepsilon_{2p_1p_2q} \}$  is a  $\mathrm{FSU}$ of  $\mathbb{K}_2$ according as  $2p_1(z+1)$   is or not a square in $\mathbb{N}$.
\item $\{\varepsilon_{2q}, \varepsilon_{p_1p_2}, \sqrt{\varepsilon_{p_1p_2}\varepsilon_{2p_1p_2q}} \}$ or  $\{\varepsilon_{2q}, \varepsilon_{p_1p_2}, \varepsilon_{2p_1p_2q} \}$  is a  $\mathrm{FSU}$ of  $\mathbb{K}_3$  according to whether $N_{\mathbb{Q}( \sqrt{p_1p_2})/\mathbb{Q}}(\varepsilon_{p_1p_2})$ equals $1$ or $-1$.
\end{enumerate}
\end{enumerate}
\end{lemma}
\begin{lemma}\label{lem3} Suppose that   $\left( \frac{2}{p_1} \right)=-\left( \frac{2}{p_2} \right)=\left( \frac{p_1}{p_2} \right)=\left( \frac{q}{p_1} \right)=-\left( \frac{q}{p_2} \right)=1$.
\begin{enumerate}[\rm1.]
\item If $2p_1(x+1)$   is a square in $\mathbb{N}$,  then
\begin{enumerate}[\rm i.]
\item $\{\varepsilon_{p_1}, \varepsilon_{2p_2q}, \sqrt{\varepsilon_{2p_1p_2q}} \}$  is a  $\mathrm{FSU}$ of  $\mathbb{K}_1$.
\item  $\{\varepsilon_{p_2}, \varepsilon_{2p_1q}, \sqrt{\varepsilon_{2p_1q}\varepsilon_{2p_1p_2q}} \}$ or $\{\varepsilon_{p_2}, \varepsilon_{2p_1q}, \varepsilon_{2p_1p_2q} \}$  is a  $\mathrm{FSU}$ of  $\mathbb{K}_2$  according as   $2p_1(z+1)$   is or not a square in $\mathbb{N}$.
\item  $\{\varepsilon_{2q}, \varepsilon_{p_1p_2}, \sqrt{\varepsilon_{p_1p_2}\varepsilon_{2p_1p_2q}} \}$ or  $\{\varepsilon_{2q}, \varepsilon_{p_1p_2}, \varepsilon_{2p_1p_2q} \}$  is a  $\mathrm{FSU}$ of  $\mathbb{K}_3$  according to whether $N_{\mathbb{Q}( \sqrt{p_1p_2} )/\mathbb{Q}}(\varepsilon_{p_1p_2})$ equals $1$ or $-1$.
\end{enumerate}
\item If $2p_2(x\pm1)$   is a square in $\mathbb{N}$,  then
\begin{enumerate}[\rm i.]
\item $\{\varepsilon_{p_1}, \varepsilon_{2p_2q}, \sqrt{\varepsilon_{2p_2q}\varepsilon_{2p_1p_2q}} \}$   is a  $\mathrm{FSU}$ of  $\mathbb{K}_1$.
\item  $\{\varepsilon_{p_2}, \varepsilon_{2p_1q}, \sqrt{\varepsilon_{2p_1p_2q}} \}$   is a  $\mathrm{FSU}$ of  $\mathbb{K}_2$.
\item  $\{\varepsilon_{2q}, \varepsilon_{p_1p_2}, \sqrt{\varepsilon_{p_1p_2}\varepsilon_{2p_1p_2q}} \}$ or  $\{\varepsilon_{q}, \varepsilon_{p_1p_2}, \varepsilon_{2p_1p_2q} \}$  is a  $\mathrm{FSU}$ of  $\mathbb{K}_3$  according to whether
    $N_{\mathbb{Q}( \sqrt{p_1p_2} )/\mathbb{Q}}(\varepsilon_{p_1p_2})$ equals $1$ or $-1$.
\end{enumerate}
{\item If $q(x\pm1)$   is a square in $\mathbb{N}$,  then
\begin{enumerate}[\rm i.]
\item  $\{\varepsilon_{p_1}, \varepsilon_{2p_2q}, \sqrt{\varepsilon_{2p_2q}\varepsilon_{2p_1p_2q}} \}$  is a  $\mathrm{FSU}$ of  $\mathbb{K}_1$.
\item  $\{\varepsilon_{p_2}, \varepsilon_{2p_1q}, \sqrt{\varepsilon_{2p_1q}\varepsilon_{2p_1p_2q}} \}$ or $\{\varepsilon_{p_2}, \varepsilon_{2p_1q}, \varepsilon_{2p_1p_2q} \}$  is a  $\mathrm{FSU}$ of  $\mathbb{K}_2$  according to whether  $q(z-1)$   is or not a square in $\mathbb{N}$.
\item  $\{\varepsilon_{2q}, \varepsilon_{p_1p_2}, \sqrt{\varepsilon_{2p_1p_2q}} \}$   is a  $\mathrm{FSU}$ of  $\mathbb{K}_3$.
\end{enumerate}}
\end{enumerate}
\end{lemma}

\begin{lemma}\label{lem4} Suppose that   $\left( \frac{2}{p_1} \right)=-\left( \frac{2}{p_2} \right)=\left( \frac{p_1}{p_2} \right)=\left( \frac{q}{p_1} \right)=\left( \frac{q}{p_2} \right)=1$.
\begin{enumerate}[\rm1.]
\item If $2p_1(x\pm1)$   is a square in $\mathbb{N}$,  then
\begin{enumerate}[\rm i.]
\item  $\{\varepsilon_{p_1}, \varepsilon_{2p_2q}, \sqrt{\varepsilon_{2p_1p_2q}} \}$  is a  $\mathrm{FSU}$ of  $\mathbb{K}_1$.
\item  $\{\varepsilon_{p_2}, \varepsilon_{2p_1q}, \sqrt{\varepsilon_{2p_1q}\varepsilon_{2p_1p_2q}} \}$ or $\{\varepsilon_{p_2}, \varepsilon_{2p_1q}, \varepsilon_{2p_1p_2q} \}$  is a  $\mathrm{FSU}$ of  $\mathbb{K}_2$  according as   $2p_1(z\pm1)$   is or not a square in $\mathbb{N}$.
\item  $\{\varepsilon_{2q}, \varepsilon_{p_1p_2}, \sqrt{\varepsilon_{p_1p_2}\varepsilon_{2p_1p_2q}} \}$ or  $\{\varepsilon_{2q}, \varepsilon_{p_1p_2}, \varepsilon_{2p_1p_2q} \}$  is a  $\mathrm{FSU}$ of  $\mathbb{K}_3$  according to whether $N_{\mathbb{Q}( \sqrt{p_1p_2} )/\mathbb{Q}}(\varepsilon_{p_1p_2})$ equals $1$ or $-1$.
\end{enumerate}
\item If $p_2(x\pm1)$   is a square in $\mathbb{N}$,  then
\begin{enumerate}[\rm i.]
\item $\{\varepsilon_{p_1}, \varepsilon_{2p_2q}, \sqrt{\varepsilon_{2p_2q}\varepsilon_{2p_1p_2q}} \}$  is a  $\mathrm{FSU}$ of  $\mathbb{K}_1$.
\item $\{\varepsilon_{p_2}, \varepsilon_{2p_1q}, \sqrt{\varepsilon_{2p_1q}\varepsilon_{2p_1p_2q}} \}$ or $\{\varepsilon_{p_2}, \varepsilon_{2p_1q}, \varepsilon_{2p_1p_2q} \}$  is a  $\mathrm{FSU}$ of  $\mathbb{K}_2$  according as   $(z\pm1)$   is or not a square in $\mathbb{N}$.
\item  $\{\varepsilon_{2q},\, \varepsilon_{p_1p_2},\, \sqrt{\varepsilon_{2q}\varepsilon_{p_1p_2}\varepsilon_{2p_1p_2q}} \}$ or  $\{\varepsilon_{2q},\, \varepsilon_{p_1p_2},\, \varepsilon_{2p_1p_2q} \}$  is a  $\mathrm{FSU}$ of  $\mathbb{K}_3$  according to whether
    $N_{\mathbb{Q}( \sqrt{p_1p_2} )/\mathbb{Q}}(\varepsilon_{p_1p_2})$ equals $1$ or $-1$.
\end{enumerate}
\item If $2q(x\pm1)$   is a square in $\mathbb{N}$,  then
\begin{enumerate}[\rm i.]
\item  $\{\varepsilon_{p_1}, \varepsilon_{2p_2q}, \sqrt{\varepsilon_{2p_2q}\varepsilon_{2p_1p_2q}} \}$  is a  $\mathrm{FSU}$ of  $\mathbb{K}_1$.
\item  $\{\varepsilon_{p_2}, \varepsilon_{2p_1q}, \sqrt{\varepsilon_{2p_1q}\varepsilon_{2p_1p_2q}} \}$ or $\{\varepsilon_{p_2}, \varepsilon_{2p_1q}, \varepsilon_{2p_1p_2q} \}$  is a  $\mathrm{FSU}$ of  $\mathbb{K}_2$  according as   $2q(z\pm1)$   is or not a square in $\mathbb{N}$.
\item  $\{\varepsilon_{2q}, \varepsilon_{p_1p_2}, \sqrt{\varepsilon_{2q}\varepsilon_{2p_1p_2q}} \}$   is a  $\mathrm{FSU}$ of  $\mathbb{K}_3$.
\end{enumerate}
\end{enumerate}
\end{lemma}

\section{The structure of the group $\mathrm{Gal}(\mathds{k}_2^{(2)}/\mathds{k})$.}
In this section we consider the field   $ \mathds{k}=\mathbb{Q}(\sqrt{2p_1p_2q} )$, where $p_1\equiv p_2\equiv -q \equiv  1 \pmod{4}$,   and the three unramified quadratic extensions\begin{center} $\mathbb{K}_1=\mathbb{Q}( \sqrt{p_1},\sqrt{2p_2q} )$, $ \mathbb{K}_2=\mathbb{Q}( \sqrt{p_2},\sqrt{2p_1q} )$ and  $ \mathbb{K}_3=\mathbb{Q}( \sqrt{2q},\sqrt{p_1p_2} )$. \end{center}
Let   $\mathrm{Cl}_2(\mathbb{K}_i)$ denote the $2$-class group of $ \mathbb{K}_i$ $(i=1, 2, 3)$.
\subsection{The metacyclic case}
\begin{theorem}\label{2}
  Let $p_1\equiv p_2\equiv-q\equiv 1 \pmod4$ be different prime integers, and $ \mathds{k}=\mathbb{Q}(\sqrt{2p_1p_2q} )$. Assume
$\mathrm{Gal}(\mathds{k}_2^{(1)}/\mathds{k})$ is a non-elementary $2$-group. Then   $\mathrm{G}=\mathrm{Gal}(\mathds{k}_2^{(2)}/\mathds{k})$ is metacyclic if and only if
$$\left( \frac{q}{p_1} \right)=\left( \frac{q}{p_2} \right)=\left( \frac{2}{p_1} \right)=\left( \frac{2}{p_2} \right)=-1.$$
More precisely,
\begin{enumerate}[\rm i.]
\item if $\left( \frac{p_1}{p_2} \right)=1$, then $\mathrm{G}$ is a metacyclic non-abelian non-modular $2$-group,
\item if $\left( \frac{p_1}{p_2} \right)=-1$, then $\mathrm{G}$ is a modular or abelian $2$-group according as  $2p_1p_2(x+1)$  is a square or not in $\mathbb{N}$.
\end{enumerate}
\end{theorem}
\begin{proof}  According to Theorems ~\ref{thm2} and  \ref{thm3}, there are five cases to distinguish. By \cite{AM} and \cite{AAM} we have:
\begin{enumerate}[\rm 1.]
\item If  $\left( \frac{2}{p_1} \right)=\left( \frac{2}{p_2} \right)=\left( \frac{p_1}{p_2} \right)=\left( \frac{q}{p_1} \right)=\left( \frac{q}{p_2} \right)=1$, then $$ \mathrm{rank}(\mathrm{Cl}_2(\mathbb{K}_1))=\mathrm{rank}(\mathrm{Cl}_2(\mathbb{K}_2)) =\mathrm{rank}(
    \mathrm{Cl}_2(\mathbb{K}_3))=3.$$
\item  If $\left( \frac{2}{p_1} \right)=\left( \frac{2}{p_2} \right)=\left( \frac{p_1}{p_2} \right)=1$ and $\left( \frac{q}{p_1} \right)=-\left( \frac{q}{p_2} \right)=1$, then $$ \mathrm{rank}(\mathrm{Cl}_2(\mathbb{K}_1))=3 \;\; \text{and} \;\; \mathrm{rank}(\mathrm{Cl}_2(\mathbb{K}_2)) =\mathrm{rank}(\mathrm{Cl}_2(\mathbb{K}_3))=2.$$
\item  If $\left( \frac{2}{p_1} \right)=\left( \frac{2}{p_2} \right)=\left( \frac{q}{p_1} \right)=\left( \frac{q}{p_2} \right)=1$ and $\left( \frac{p_1}{p_2} \right)=-1$, then $$ \mathrm{rank}(\mathrm{Cl}_2(\mathbb{K}_1))= \mathrm{rank}(\mathrm{Cl}_2(\mathbb{K}_2)) =2  \;\; \text{and} \;\; \mathrm{rank}(\mathrm{Cl}_2(\mathbb{K}_3))=3.$$
\item  If $\left( \frac{2}{p_1} \right)=-\left( \frac{2}{p_2} \right)=1$ and $\left( \frac{p_1}{p_2} \right)=\left( \frac{q}{p_1} \right)=1$, then
$$ \mathrm{rank}(\mathrm{Cl}_2(\mathbb{K}_1))=3 \;\; \text{and} \;\; \mathrm{rank}(\mathrm{Cl}_2(\mathbb{K}_2)) =\mathrm{rank}(\mathrm{Cl}_2(\mathbb{K}_3))=2.$$
\item  If $\left( \frac{2}{p_1} \right)=\left( \frac{2}{p_2} \right)=\left( \frac{q}{p_1} \right)=\left( \frac{q}{p_2} \right)=-1$, then
\begin{enumerate}[\rm i.]
\item  If $\left(\frac{p_1}{p_2}\right)=1$, then $$ \mathrm{rank}(\mathrm{Cl}_2(\mathbb{K}_1))=\mathrm{rank}(\mathrm{Cl}_2(\mathbb{K}_2)) =\mathrm{rank}(\mathrm{Cl}_2(\mathbb{K}_3))=2.$$
\item If $\left(\frac{p_1}{p_2}\right)=-1$, then  $$ \mathrm{rank}(\mathrm{Cl}_2(\mathbb{K}_1))= \mathrm{rank}(\mathrm{Cl}_2(\mathbb{K}_2)) =1  \; \text{and} \; \mathrm{rank}(\mathrm{Cl}_2(\mathbb{K}_3))=2.$$
\end{enumerate}
\end{enumerate}
Hence the results are   deduced from Theorem ~\ref{thm0}, Lemma ~\ref{lem7} and \cite[ Theorem 2]{B.L.S-98}.
\end{proof}
\subsection{The non-metacyclic case}

  Assuming $\mathrm{Cl}_2(\mathds{k})\simeq (2, 2^n)$, $n\geq2$,
for the non-metacyclic case we have four cases to distinguish, according to Theorems \ref{thm3} and  \ref{2}. For simplicity, we will  denote by  $q_i$ the unit index of the field $\mathbb{K}_i$ $(i=1, 2, 3)$. In all that follows, we use the notations of \cite[page 336]{Kaplan}. Put  $p_1=2e^2+(-1)^{\gamma} d^2$, $q=2r^2+(-1)^{\gamma} s^2$ and $A=sd+2er+2\gamma(es+dr)$ according as $\left( \frac{2}{q} \right)=(-1)^{\gamma+1}$, where $\gamma\in\{0, 1\}$.

\subsubsection{\bf Case $1$: $\left( \frac{2}{p_1} \right)=\left( \frac{2}{p_2} \right)=\left( \frac{p_1}{p_2} \right)=\left( \frac{q}{p_1} \right)=-\left( \frac{q}{p_2} \right)=1$ }
\begin{theorem}\label{thm8} Let $\delta \in \{1, p_1, 2p_1\}$ be such that $\delta(x\pm1)$   is a square in  $\mathbb{N}$.   The group  $\mathrm{Cl}_2(\mathds{k}^{(1)}_2)\simeq\mathrm{Gal}(\mathds{k}_2^{(2)}/\mathds{k}_2^{(1)})$ is non-elementary cyclic if and only if one of the two
following assertions holds:
\begin{enumerate}[\rm I.]
\item
\begin{enumerate}[\rm i.]
\item   $\delta(z\pm1)$ is not a square in  $\mathbb{N}$,
\item   at least one of the elements  $\left\{\left(\frac{A}{p_1}\right) ,\left(\frac{2q}{p_1}\right)_4\right\}$  equals  $-1$, and
\item     either 
\begin{enumerate}[\rm a.]
    \item $\delta=1$ and $\left(\frac{p_1}{p_2}\right)_4 = \left(\frac{p_2}{p_1}\right)_4$, or
    \item $\delta\neq 1$ and $\left(\frac{p_1}{p_2}\right)_4 = \left(\frac{p_2}{p_1}\right)_4 =1$.
\end{enumerate}
\end{enumerate}
\item
\begin{enumerate}[\rm i.]
\item   $\delta(z\pm1)$   is a square in  $\mathbb{N}$ or $\left(\frac{A}{p_1}\right)=\left(\frac{2q}{p_1}\right)_4=1$, and
\item  either 
\begin{enumerate}[\rm a.]
    \item $\delta=1$ and $\left(\frac{p_1}{p_2}\right)_4 \neq \left(\frac{p_2}{p_1}\right)_4$, or
    \item $\delta\neq 1$ and one of $\left(\frac{p_1}{p_2}\right)_4, \left(\frac{p_2}{p_1}\right)_4$ is equal to $-1$.
\end{enumerate}

\end{enumerate}
\end{enumerate}
\end{theorem}
\begin{proof}Form Theorem ~\ref{thm4}, we must calculate the $2$-class numbers of $\mathbb{K}_i$.
\begin{itemize}
    \item By \cite{Kaplan}, if $\left( \frac{q}{p_2} \right)=-1$, then $h_2(2p_2q)=2$, and,  according to Lemma ~\ref{lem2}, $q_1=2$. In this case, the  $2$-class number of $ \mathbb{K}_1$ is given by \cite{Wada}:
$$ h_2(\mathbb{K}_1)=\frac{1}{4}  q_1 h_2(p_1) h_2(2p_2q)h_2(2p_1p_2q)=h_2(2p_1p_2q).$$
    \item If $\left( \frac{2}{p_1} \right)=\left( \frac{q}{p_1} \right)=1$, then, by \cite{Kaplan}, $h_2(2p_1q)\geq 4$. More precisely, $h_2(2p_1q)=4$
  if and only if  at least one of the elements  $\left(\frac{A}{p_1}\right) ,\left(\frac{2q}{p_1}\right)_4$  equals  $-1$.
    The $2$-class number of  $ \mathbb{K}_2$ is given by:
$$ h_2(\mathbb{K}_2)=\frac{1}{4} q_2 h_2(p_2) h_2(2p_1q)h_2(2p_1p_2q)=\frac{1}{4} q_2 h_2(2p_1q)h_2(2p_1p_2q),$$
so $ h_2(\mathbb{K}_2)=h_2(2p_1p_2q)$  if and only if  $ q_2=1$ and $h_2(2p_1q)=4$. On the other hand, by Lemma ~\ref{lem2}, $ q_2=1$   if and only if  $\delta(z\pm1)$ is not a square in  $\mathbb{N}$.
    \item Similarly, the $2$-class number of  $ \mathbb{K}_3$  is given by:
$$ h_2(\mathbb{K}_3)=\frac{1}{4} q_3 h_2(2q) h_2(p_1p_2)h_2(2p_1p_2q)=\frac{1}{4} q_3 h_2(p_1p_2)h_2(2p_1p_2q),$$
so $ h_2(\mathbb{K}_3)=h_2(2p_1p_2q)$  if and only if  either $q_3=1$ and $h_2(p_1p_2)=4$, or $q_3=2$ and $h_2(p_1p_2)=2$. In this case, according to \cite{Kucera} (see also \cite{Sc34}) and Lemma ~\ref{lem2}, $ h_2(\mathbb{K}_3)=h_2(2p_1p_2q)$  if and only if one of the two following conditions is satisfied:
\begin{enumerate}[\rm1.]
\item   $\delta=1$ and $\left(\frac{p_1}{p_2}\right)_4 \neq \left(\frac{p_2}{p_1}\right)_4$.
\item   $\delta\neq 1$ and either $ \left(\frac{p_1}{p_2}\right)_4=-1$ or $\left(\frac{p_2}{p_1}\right)_4=-1$.
\end{enumerate}
\end{itemize}
Using Theorem ~\ref{thm4}, we get the results.
\end{proof}
\begin{example}\label{9}
Put $\alpha=\left(\frac{A}{p_1}\right)$, $s=\left(\frac{2q}{p_1}\right)_4$, $t_1=\left(\frac{p_1}{p_2}\right)_4$, $t_2=\left(\frac{p_2}{p_1}\right)_4$, $c=\mathrm{Cl}(\mathds{k}_2^{(1)})$, $n=h_2(\mathds{k})$, $n_i=h_2(\mathbb{K}_i)$ $(i=1, 2, 3)$ and $q_0=q(\mathbb{K}_0/\mathds{k})$, and by using PARI/GP \cite{GP-16}, we get the following examples for the case: $x\pm1$ is a square, $z+1$ and $z-1$ are not squares, $(\alpha=-1$ or $s=-1)$ and $t_1=t_2$.
\begin{longtable}{| c | c | c | c | c |c | c | c | c| c | c | c |c|c|c|c|}
\hline
$d=2p_1p_2q  $&$ q_1 $&$ q_2 $&$ q_3 $&$ \alpha $&$ s $&$ t_1 $&$ t_2  $&$ n $&$ n_1 $&$ n_2 $&$ n_3 $&$ q_0$&$ c$\\
\hline
\endfirsthead
\hline
$d=2p_1.p_2.q  $&$ q_1 $&$ q_2 $&$ q_3 $&$ \alpha $&$ s $&$ t_1 $&$ t_2  $&$ n $&$ n_1 $&$ n_2 $&$ n_3 $&$ q_0$&$ c$\\
\hline
\endhead
$38982=2\cdot 73\cdot 89\cdot 3$&$2$&$1$&$2$&$-1$&$-1$&$1$&$1$&$8$&$8$&$8$&$16$&$16$&$[4]$\\
$60006=2\cdot 73\cdot 137\cdot 3$&$2$&$1$&$2$&$-1$&$-1$&$1$&$1$&$8$&$8$&$8$&$64$&$16$&$[16]$\\
$298862=2\cdot 73\cdot 89\cdot 23$&$2$&$1$&$2$&$-1$&$-1$&$1$&$1$&$8$&$8$&$8$&$16$&$16$&$[12]$\\
\hline
\end{longtable}

\end{example}
\begin{theorem}\label{thm9}  Let $\delta \in \{1, p_1, 2p_1\}$  such that  $\delta(x\pm1)$   is a square in  $\mathbb{N}$.
 Then $\#\mathrm{Cl}_2(\mathds{k}^{(1)}_2)=2$  if and only if   the following conditions are satisfied:

\begin{enumerate}[\rm i.]
\item   $\delta(z\pm1)$ is not a square in  $\mathbb{N}$,
\item  at least one of the elements  $\left\{\left(\frac{A}{p_1}\right) ,\left(\frac{2q}{p_1}\right)_4\right\}$  equals  $-1$.
\item $ \left(\frac{p_1}{p_2}\right)_4 \neq \left(\frac{p_2}{p_1}\right)_4 $.
\end{enumerate}
\end{theorem}

\begin{proof} Suppose that  $\#\mathrm{Cl}_2(\mathds{k}^{(1)}_2)=2$. Then, according to Corollary  ~\ref{coro5}, $h_2(\mathbb{K}_i)=h_2(\mathds{k})$ for all $i=1, 2, 3$. By the proof of Theorem ~\ref{thm8}, the equality  $h_2(\mathbb{K}_2)=h_2(\mathds{k})$ implies the two first conditions and $q_2=1$. On the other hand, as $ q_1=2$, from Theorem ~\ref{thm7} we infer $q_3=2$. Accordingly, $h_2(p_1p_2)=2$ and $N_{\mathbb{Q}( \sqrt{p_1p_2} )/\mathbb{Q}}(\varepsilon_{p_1p_2})=1$, which is equivalent to  $ \left(\frac{p_1}{p_2}\right)_4 \neq \left(\frac{p_2}{p_1}\right)_4 $ (see \cite{Sc34}).

Reciprocally, suppose the three   conditions (i), (ii) and (iii) are satisfied. Applying results of the proof of Theorem  ~\ref{thm8}, we get $h_2(\mathbb{K}_i)=h_2(\mathds{k})$ for all $i=1, 2, 3$.  Let $\mathbb{K}_0=\mathbb{K}_1\mathbb{K}_2\mathbb{K}_3=\mathbb{Q}(\sqrt{p_1},\sqrt{p_2},\sqrt{2q})$,  and denote by $E_{\mathbb{K}_i}$ the unit group of $\mathbb{K}_i$ and by  $q(\mathbb{K}_0/\mathds{k})=[E_{\mathbb{K}_0}:E_{\mathbb{K}_1}E_{\mathbb{K}_2}E_{\mathbb{K}_3}]$  the unit index of  $ \mathbb{K}_0/\mathds{k}$. Hence, by Corollary  ~\ref{coro5}, it remains to prove only $q(\mathbb{K}_0/\mathds{k})=4$. According to Lemma ~\ref{lem2}, we have:
\begin{enumerate}[\rm1.]
\item $E_{\mathbb{K}_1}=\langle -1,\varepsilon_{p_1}, \varepsilon_{2p_2q},\varepsilon\rangle$, where $\varepsilon=\sqrt{\varepsilon_{2p_1p_2q}}$ or $\sqrt{\varepsilon_{2p_2q}\varepsilon_{2p_1p_2q}}$  according to  $2p_1(x\pm1)$   is or not a square in $\mathbb{N}$,
\item $E_{\mathbb{K}_2}=\langle-1,\varepsilon_{p_2}, \varepsilon_{2p_1q},\varepsilon_{2p_1p_2q}\rangle$,
\item $E_{\mathbb{K}_3}=\langle-1,\varepsilon_{2q}, \varepsilon_{p_1p_2},\varepsilon'\rangle$, where $\varepsilon'=\sqrt{\varepsilon_{2q}\varepsilon_{2p_1p_2q}}$,  $\sqrt{\varepsilon_{2q}\varepsilon_{p_1p_2}\varepsilon_{2p_1p_2q}}$ or $\sqrt{\varepsilon_{p_1p_2}\varepsilon_{2p_1p_2q}}$  according to  $(x\pm1)$,  $p_1(x\pm1)$ or $2p_1(x\pm1)$  is a square in  $\mathbb{N}$.
\end{enumerate}
 Put

$$\begin{array}{ll}
A&=\varepsilon_{p_1}^{a_1}\varepsilon_{2p_2q}^{a_2}\varepsilon^{a_3}, \;\;a_1,a_2,a_3 \in \{0,1\},\\\\
B&=\varepsilon_{p_2}^{b_1}\varepsilon_{2p_1q}^{b_2}\varepsilon_{2p_1p_2q}^{b_3},\;\; b_1,b_2,b_3 \in \{0,1\},\\\\
C&=\varepsilon_{2q}^{c_1}\varepsilon_{p_1p_2}^{c_2}\varepsilon'^{c_3},\;\; c_1,c_2,c_3 \in \{0,1\},\\\\
\eta^2&=\pm A.B.C.
\end{array}$$

So
$$\begin{array}{ll}
N_{\mathbb{K}_0/\mathbb{K}_1}(\eta^2)&=(-1)^{b_1}(\pm \varepsilon_{2p_1p_2q})^{c_3}(\varepsilon_{2p_1p_2q}^{b_3} A)^2,\\\\
N_{\mathbb{K}_0/\mathbb{K}_2}(\eta^2)&=(-1)^{a_1}(\pm \varepsilon_{2p_1p_2q})^{a_3}(\pm \varepsilon_{2p_1p_2q})^{c_3} B^2,\\\\
N_{\mathbb{K}_0/\mathbb{K}_3}(\eta^2)&=(-1)^{a_1}(-1)^{b_1}(\pm \varepsilon_{2p_1p_2q})^{a_3}(\varepsilon_{2p_1p_2q}^{b_3} C)^2.

\end{array}$$

Assume $\eta \in \mathbb{K}_0$, if $a_3\neq0$ or $c_3\neq0$, then $\sqrt{\varepsilon_{2p_1p_2q}} \in \mathbb{K}_2$ or $\sqrt{\varepsilon_{2p_1p_2q}} \in \mathbb{K}_3$, which contradicts  Lemma ~\ref{lem2}. On the other hand, if $a_3=c_3=0$ and ($a_1=1$ or $b_1=1$ ), then $N_{\mathbb{K}_0/\mathbb{K}_1}(\eta^2)<0$ or $N_{\mathbb{K}_0/\mathbb{K}_2}(\eta^2)<0$, which contradicts the fact that $N_{\mathbb{K}_0/\mathbb{K}_i}(\eta^2)>0$. Therefore, $a_1=b_1=a_3=c_3=0$ and we get  $$\eta^2=\pm\varepsilon_{2p_2q}^{a_2}\varepsilon_{2p_1q}^{b_2}\varepsilon_{2p_1p_2q}^{b_3}\varepsilon_{2q}^{c_1}\varepsilon_{p_1p_2}^{c_2}.$$
From the proof of Lemma ~\ref{lem2}, we deduce thet $\sqrt{\varepsilon_{2q}\varepsilon_{2p_iq}}$ or $\sqrt{\varepsilon_{2p_iq}} \in E_{\mathbb{K}_0}$ $(i=1,2)$. According to our  assumption,  $N_{\mathbb{Q}( \sqrt{p_1p_2} )/\mathbb{Q}}(\varepsilon_{p_1p_2})=1$, which implies that $\sqrt{\varepsilon_{p_1p_2}} \in E_{\mathbb{K}_0}$.
We distinguish the following cases:
\begin{enumerate}[\rm i.]
\item If  $(x\pm1)$  is a square in  $\mathbb{N}$, then $\sqrt{\varepsilon_{2q}},\sqrt{\varepsilon_{2p_2q}},\sqrt{\varepsilon_{2p_1p_2q}} \notin E_{\mathbb{K}_0} $,  and Lemma ~\ref{lem2} implies that 
\begin{align*}
    \sqrt{\varepsilon_{p_1p_2}} &\notin E_{\mathbb{K}_1}E_{\mathbb{K}_2}E_{\mathbb{K}_3}, \\
    \sqrt{\varepsilon_{2q}\varepsilon_{2p_1q}},\sqrt{\varepsilon_{2p_1q}} &\notin E_{\mathbb{K}_1}E_{\mathbb{K}_2}E_{\mathbb{K}_3}, \\
    \sqrt{\varepsilon_{2q}\varepsilon_{2p_1p_2q}},\sqrt{\varepsilon_{2p_2q}\varepsilon_{2p_1p_2q}}&\in E_{\mathbb{K}_1}E_{\mathbb{K}_2}E_{\mathbb{K}_3}.
\end{align*} By a case by case study, we obtain
$$E_{\mathbb{K}_0}=\langle \sqrt{\varepsilon_{2p_1q}}\;\; \text{or} \;\; \sqrt{\varepsilon_{2q}\varepsilon_{2p_1q}}, \sqrt{\varepsilon_{p_1p_2}},E_{\mathbb{K}_1}E_{\mathbb{K}_2}E_{\mathbb{K}_3} \rangle.$$
\item If  $p_1(x\pm1)$  is a square in  $\mathbb{N}$, then $\sqrt{\varepsilon_{2q}},\sqrt{\varepsilon_{2p_2q}},\sqrt{\varepsilon_{2p_1p_2q}} \notin E_{\mathbb{K}_0} $,  By Lemma ~\ref{lem2}, we have 
\begin{align*}
    \sqrt{\varepsilon_{p_1p_2}} &\notin E_{\mathbb{K}_1}E_{\mathbb{K}_2}E_{\mathbb{K}_3}, \\
    \sqrt{\varepsilon_{2q}\varepsilon_{2p_1q}},\sqrt{\varepsilon_{2p_1q}} &\notin E_{\mathbb{K}_1}E_{\mathbb{K}_2}E_{\mathbb{K}_3}, \\
    \sqrt{\varepsilon_{2q}\varepsilon_{p_1p_2}\varepsilon_{2p_1p_2q}},\sqrt{\varepsilon_{2p_2q}\varepsilon_{2p_1p_2q}} &\in E_{\mathbb{K}_1}E_{\mathbb{K}_2}E_{\mathbb{K}_3}.
\end{align*}

Thus we obtain
$$E_{\mathbb{K}_0}=\langle\sqrt{\varepsilon_{2p_1q}}\;\; \text{or} \;\; \sqrt{\varepsilon_{2q}\varepsilon_{2p_1q}}, \sqrt{\varepsilon_{p_1p_2}},E_{\mathbb{K}_1}E_{\mathbb{K}_2}E_{\mathbb{K}_3}\rangle.$$
\item If   $2p_1(x\pm1)$  is a square in   $\mathbb{N}$, then $\sqrt{\varepsilon_{2q}},\sqrt{\varepsilon_{2p_2q}}\notin E_{\mathbb{K}_0} $ and $\sqrt{\varepsilon_{2q}\varepsilon_{2p_2q}}\in E_{\mathbb{K}_0} $, by Lemma ~\ref{lem2}, $\sqrt{\varepsilon_{p_1p_2}},\sqrt{\varepsilon_{2p_1p_2q}} \in E_{\mathbb{K}_1}E_{\mathbb{K}_2}E_{\mathbb{K}_3} $, from which we deduce that
$$E_{\mathbb{K}_0}=\langle\sqrt{\varepsilon_{2p_1q}}\;\; \text{or} \;\; \sqrt{\varepsilon_{2q}\varepsilon_{2p_1q}}, \sqrt{\varepsilon_{2q}\varepsilon_{2p_2q}},E_{\mathbb{K}_1}E_{\mathbb{K}_2}E_{\mathbb{K}_3}\rangle.$$
\end{enumerate}
In the three cases we get $q(\mathbb{K}_0/\mathds{k})=4$, so it suffices to apply Corollary ~\ref{coro5} to obtain the results.
\end{proof}
\begin{example} Keep the notation of Example \ref{9}. For the case  $p_1(x\pm1)$ is a square, $p_1(z\pm1)$ is not a square, $(\alpha=-1$ or $s=-1)$ and $t_1\neq t_2$, we have
\begin{longtable}{| l | c | c | c | c |c | c | c | c| c | c | c | c|c|c|c|}
\hline
$d=p_1p_2q  $&$ q_1 $&$ q_2 $&$ q_3 $&$ \alpha $&$ s $&$ t_1 $&$ t_2  $&$ n $&$ n_1 $&$ n_2 $&$ n_3 $&$ q_0 $&$ c$\\
\hline
\endfirsthead
\hline
$d=p_1.p_2.q  $&$ q_1 $&$ q_2 $&$ q_3 $&$ \alpha $&$ s $&$ t_1 $&$ t_2  $&$ n $&$ n_1 $&$ n_2 $&$ n_3 $&$ q_0 $&$ c$\\
\hline
\endhead
$51798=2\cdot 97\cdot 89\cdot 3$&$2$&$1$&$2$&$-1$&$1$&$1$&$-1$&$8$&$8$&$8$&$8$&$16$&$[2]$\\
$64862=2\cdot 113\cdot 41\cdot 7$&$2$&$1$&$2$&$-1$&$1$&$1$&$-1$&$8$&$8$&$8$&$8$&$16$&$[6]$\\
$113734=2\cdot 73\cdot 41\cdot 19$&$2$&$1$&$2$&$1$&$-1$&$-1$&$1$&$8$&$8$&$8$&$8$&$16$&$[6]$\\
\hline
\end{longtable}

\end{example}
\subsubsection{\bf Case $2$: $\left( \frac{2}{p_1} \right)=\left( \frac{2}{p_2} \right)=\left( \frac{q}{p_1} \right)=\left( \frac{q}{p_2} \right)=-\left( \frac{p_1}{p_2} \right)=1$ }

\begin{lemma} \label{lem8} Let  $p_1 \equiv p_2\equiv-q \equiv1 \pmod4$ be three positive prime integers satisfying  $\left( \frac{2}{p_1} \right)=\left( \frac{2}{p_2} \right)=\left( \frac{q}{p_1} \right)=\left( \frac{q}{p_2} \right)=-\left( \frac{p_1}{p_2} \right)=1$. Then the rank of the $2$-class group  of $ \mathbb{K}_0=\mathbb{Q}( \sqrt{p_1},\sqrt{p_2},\sqrt{2q} )$  equals  $3$.
\end{lemma}
\begin{proof} As $\left( \frac{p_1}{p_2} \right)=-1$, it is well known (cf.\ \cite{Kucera}) that the class number of  $\mathds{k}_0=\mathbb{Q}( \sqrt{p_1},\sqrt{p_2} )$ is odd. Consider the extension $ \mathbb{K}_0/\mathds{k}_0$. Then according to \cite{Gr}, the rank of the $2$-class group  of  $\mathbb{K}_0$  is given by the formula:
$$\text{rank}(\mathrm{Cl}_2(\mathbb{K}_0))=r-e-1 ,$$
where $r$ is the number of finite and infinite primes of $\mathds{k}_0$ that ramify in  $\mathbb{K}_0/\mathds{k}_0$ and $e$ is defined by $2^{e}=[E_{\mathds{k}_0}:E_{\mathds{k}_0} \cap N_{\mathbb{K}_0/\mathds{k}_0}(\mathbb{K}_0^*)]\leq 2^4$. As  $r=8$ ($4$ primes above   $2$ and $4$ above $q$), $\text{rank}(\mathrm{Cl}_2(\mathbb{K}_0))=r-e-1=8-e-1\geq3$. On the other hand, the Schreier's inequality implies that $\text{rank}(\mathrm{Cl}_2(\mathbb{K}_0))\leq 3$, concluding the proof.
\end{proof}
\begin{theorem} \label{thm10} Let  $p_1 \equiv p_2\equiv-q \equiv1 \pmod4$ be three positive prime integers satisfying  $\left( \frac{2}{p_1} \right)=\left( \frac{2}{p_2} \right)=\left( \frac{q}{p_1} \right)=\left( \frac{q}{p_2} \right)=-\left( \frac{p_1}{p_2} \right)=1$. Then
$$\mathrm{rank}( \mathrm{Cl}_2(\mathds{k}^{(1)}_2))\geq 2.$$
\end{theorem}
\begin{proof} Suppose that  $\mathrm{rank}( \mathrm{Cl}_2(\mathds{k}^{(1)}_2))=1$, then, by Lemma ~\ref{lem6}, the proof of Theorem ~\ref{2} and class field theory, the Galois  groups $\mathrm{Gal}(\mathds{k}_2^{(2)}/\mathbb{K}_1)$ and $\mathrm{Gal}(\mathds{k}_2^{(2)}/\mathbb{K}_2)$ are metacyclic, and since  $\mathbb{K}_0$ is an unramified quadratic  extension of  $\mathbb{K}_1$, then Lemma ~\ref{lem7} implies that $\text{rank}(\mathrm{Cl}_2(\mathbb{K}_0))\leq 2$, which contradicts  Lemma ~\ref{lem8}.
\end{proof}
\begin{example}
 For $d=47158=2\cdot 73\cdot 17\cdot 19$, we have  $\mathrm{Cl}_2(\mathds{k}^{(1)}_2)$ is of type $(2, 4)$, and for
$d=59942=2\cdot 17\cdot 41\cdot 43$, we have $\mathrm{Cl}_2(\mathds{k}^{(1)}_2)$ is of type $(2, 4)$.
\end{example}
\subsubsection{\bf Case $3$: $\left( \frac{2}{p_1} \right)=-\left( \frac{2}{p_2} \right)=\left( \frac{p_1}{p_2} \right)=\left( \frac{q}{p_1} \right)=-\left( \frac{q}{p_2} \right)=1$ }

Using similar arguments as above, we prove the following two theorems.
\begin{theorem}\label{thm11}  Let $\delta \in \{2p_1, 2p_2, q\}$ be such that  $\delta(x\pm1)$   is a square in  $\mathbb{N}$. The group $\mathrm{Cl}_2(\mathds{k}^{(1)}_2)\simeq\mathrm{Gal}(\mathds{k}_2^{(2)}/\mathds{k}_2^{(1)})$ is cyclic non-elementary   if and only if   one of the following two conditions is satisfied:
\begin{enumerate}[\rm I.]
\item
\begin{enumerate}[\rm i.]
\item  $\delta \neq 2p_2$ and $\delta(z\pm1)$ is not a square in  $\mathbb{N}$, and
\item  at least one of the elements  $\left\{\left(\frac{A}{p_1}\right) ,\left(\frac{2q}{p_1}\right)_4\right\}$  equals  $-1$, and
\item     either
\begin{enumerate}[\rm a.]
    \item $\delta=q$ and $\left(\frac{p_1}{p_2}\right)_4 = \left(\frac{p_2}{p_1}\right)_4$, or
    \item $\delta\neq q $ and $ \left(\frac{p_1}{p_2}\right)_4 = \left(\frac{p_2}{p_1}\right)_4 =1$
\end{enumerate}
\end{enumerate}
\item
\begin{enumerate}[\rm i.]
\item  $\delta=2p_2$ or $\delta(z\pm1)$   is a square in  $\mathbb{N}$ or $\left(\frac{A}{p_1}\right)=\left(\frac{2q}{p_1}\right)_4=1$, and
\item either
\begin{enumerate}[\rm a.]
    \item $\delta=q$ and $\left(\frac{p_1}{p_2}\right)_4 \neq \left(\frac{p_2}{p_1}\right)_4$, or
    \item $\delta\neq q $ and  one of $ \left(\frac{p_1}{p_2}\right)_4, \left(\frac{p_2}{p_1}\right)_4 $ is equal to $-1$.
\end{enumerate}

\end{enumerate}
\end{enumerate}
\end{theorem}

\begin{example} {Keep the notation of Example \ref{9}. For the case $2p_2(x\pm1)$ is a square
and $t_1=-1$ or $t_2=-1$, we have}\nopagebreak

\setlength\tabcolsep{5pt}
\begin{tabular}{| c | c | c | c | c |c | c | c | c| c | c | c |c|c|c|c|}
\hline
$d=2p_1p_2q  $&$ q_1 $&$ q_2 $&$ q_3 $&$ \alpha $&$ s $&$ t_1 $&$ t_2  $&$ n $&$ n_1 $&$ n_2 $&$ n_3 $&$ q_0$&$ c$\\
\hline
$17630=2\cdot 41\cdot 5\cdot 43$&$2$&$2$&$2$&$1$&$1$&$1$&$-1$&$8$&$8$&$32$&$8$&$32$&$[8]$\\
$29614=2\cdot 17\cdot 13\cdot 67$&$2$&$2$&$2$&$1$&$-1$&$-1$&$1$&$8$&$8$&$16$&$8$&$32$&$[4]$\\
$34238=2\cdot 17\cdot 53\cdot 19$&$2$&$2$&$1$&$1$&$1$&$-1$&$-1$&$8$&$8$&$32$&$8$&$16$&$[8]$\\
$41830=2\cdot 89\cdot 5\cdot 47$&$2$&$2$&$1$&$-1$&$-1$&$-1$&$-1$&$16$&$16$&$32$&$16$&$16$&$[4]$\\
$59630=2\cdot 89\cdot 5\cdot 67$&$2$&$2$&$1$&$-1$&$1$&$-1$&$-1$&$8$&$8$&$16$&$8$&$16$&$[4]$\\
$69782=2\cdot 41\cdot 37\cdot 23$&$2$&$2$&$2$&$1$&$-1$&$-1$&$1$&$8$&$8$&$16$&$8$&$32$&$[4]$\\
$91078=2\cdot 113\cdot 13\cdot 31$&$2$&$2$&$2$&$-1$&$-1$&$1$&$-1$&$16$&$16$&$32$&$16$&$32$&$[12]$\\
\hline
\end{tabular}

\end{example}

\begin{theorem}\label{thm12}  Let $\delta \in \{2p_1, 2p_2, q\}$ be such that  $\delta(x\pm1)$   is a square in  $\mathbb{N}$. The order $\#\mathrm{Cl}_2(\mathds{k}^{(1)}_2)=2$  if and only if   the following conditions are satisfied:

\begin{enumerate}[\rm i.]
\item  $\delta \neq 2p_2$ and $\delta(z\pm1)$ is not a square in  $\mathbb{N}$,
\item  at least one of the elements  $\left\{\left(\frac{A}{p_1}\right) ,\left(\frac{2q}{p_1}\right)_4\right\}$  equals  $-1$,
\item $ \left(\frac{p_1}{p_2}\right)_4 \neq \left(\frac{p_2}{p_1}\right)_4 $.
\end{enumerate}
\end{theorem}

\begin{example} Keep the notation of Example \ref{9}. For the case $q(x\pm1)$ is a square, $q(z+1)$ and $q(z-1)$ are not squares, $(\alpha=-1$ or $s=-1)$ and $t_1\neq t_2$ we have
\begin{longtable}{| c | c | c | c | c |c | c | c | c| c | c | c |c|c|c|c|}
\hline
$d=2p_1p_2q  $&$ q_1 $&$ q_2 $&$ q_3 $&$ \alpha $&$ s $&$ t_1 $&$ t_2  $&$ n $&$ n_1 $&$ n_2 $&$ n_3 $&$ q_0$&$ c$\\
\hline
\endfirsthead
\hline
$d=2p_1.p_2.q  $&$ q_1 $&$ q_2 $&$ q_3 $&$ \alpha $&$ s $&$ t_1 $&$ t_2  $&$ n $&$ n_1 $&$ n_2 $&$ n_3 $&$ q_0$&$ c$\\
\hline
\endhead
$9430=2\cdot 41\cdot 5\cdot 23$&$2$&$1$&$2$&$1$&$-1$&$1$&$-1$&$16$&$16$&$16$&$16$&$16$&$[2]$\\
$20774=2\cdot 17\cdot 13\cdot 47$&$2$&$1$&$2$&$1$&$-1$&$-1$&$1$&$8$&$8$&$8$&$8$&$16$&$[2]$\\
$94054=2\cdot 41\cdot 37\cdot 31$&$2$&$1$&$2$&$1$&$-1$&$-1$&$1$&$8$&$8$&$8$&$8$&$16$&$[2]$\\
$102638=2\cdot 73\cdot 37\cdot 19$&$2$&$1$&$2$&$1$&$-1$&$-1$&$1$&$8$&$8$&$8$&$8$&$16$&$[6]$\\
\hline
\end{longtable}
\end{example}
\subsubsection{\bf Case $4$: $\left( \frac{2}{p_1} \right)=-\left( \frac{2}{p_2} \right)=\left( \frac{p_1}{p_2} \right)=\left( \frac{q}{p_1} \right)=\left( \frac{q}{p_2} \right)=1$}

Using similar arguments as above, we prove the following two theorems.
\begin{theorem}\label{thm13}  Let $\delta \in \{2p_1, p_2, 2q\}$  such that  $\delta(x\pm1)$   is a square in  $\mathbb{N}$. The group $\mathrm{Cl}_2(\mathds{k}^{(1)}_2)\simeq\mathrm{Gal}(\mathds{k}_2^{(2)}/\mathds{k}_2^{(1)})$ is cyclic non-elementary   if and only if   one of the following two conditions is satisfied:
\begin{enumerate}[\rm I.]
\item
\begin{enumerate}[\rm i.]
\item $ \delta=p_2\ ($resp. $\delta\neq p_2)$ and $(z\pm1)$  $($resp. $\delta(z\pm1))$ is not a square in  $\mathbb{N}$,
\item  at least one of the elements  $\left\{\left(\frac{A}{p_1}\right) ,\left(\frac{2q}{p_1}\right)_4\right\}$  equals  $-1$,
\item     either 
\begin{enumerate}[\rm a.]
    \item $\delta=2q$ and $\left(\frac{p_1}{p_2}\right)_4 = \left(\frac{p_2}{p_1}\right)_4$, or
    \item $\delta\neq 2q$ and $\left(\frac{p_1}{p_2}\right)_4= \left(\frac{p_2}{p_1}\right)_4=1$.
\end{enumerate}    
\end{enumerate}
\item
\begin{enumerate}[\rm i.]
\item  $(z\pm1)\  ($resp. $\delta(z\pm1))$ is a square in  $\mathbb{N}$ if $ \delta=p_2\ ($resp. $\delta\neq p_2)$ or $\left(\frac{A}{p_1}\right)=\left(\frac{2q}{p_1}\right)_4=1$,
\item  either 
\begin{enumerate}[\rm a.]
    \item $\delta=2q$ and $\left(\frac{p_1}{p_2}\right)_4 \neq \left(\frac{p_2}{p_1}\right)_4$, or
    \item $\delta\neq 2q$ and one of $\left(\frac{p_1}{p_2}\right)_4, \left(\frac{p_2}{p_1}\right)_4$ is equal to $-1$.
\end{enumerate} 
\end{enumerate}
\end{enumerate}
\end{theorem}
\begin{example}Keep the notation of Example \ref{9}. For the case $p_2(x\pm1)$ is a square, ($z\pm1$ is a square or $(\alpha=s=1))$ and $(t_1=-1$ or  $t_2=-1)$, we have
\begin{longtable}{| c | c | c | c | c |c | c | c | c| c | c | c |c|c|c|c|}
\hline
$d=2p_1p_2q  $&$ q_1 $&$ q_2 $&$ q_3 $&$ \alpha $&$ s $&$ t_1 $&$ t_2  $&$ n $&$ n_1 $&$ n_2 $&$ n_3 $&$ q_0$&$ c$\\
\hline
\endfirsthead
\hline
$d=2p_1.p_2.q  $&$ q_1 $&$ q_2 $&$ q_3 $&$ \alpha $&$ s $&$ t_1 $&$ t_2  $&$ n $&$ n_1 $&$ n_2 $&$ n_3 $&$ q_0$&$ c$\\
\hline
\endhead
$84422=2\cdot 17\cdot 13\cdot 191$&$2$&$2$&$2$&$1$&$-1$&$-1$&$1$&$8$&$8$&$16$&$8$&$32$&$[12]$\\
$113102=2\cdot 97\cdot 53\cdot 11$&$2$&$1$&$2$&$1$&$1$&$1$&$-1$&$8$&$8$&$16$&$8$&$16$&$[4]$\\
$123710=2\cdot 89\cdot 5\cdot 139$&$2$&$1$&$1$&$1$&$1$&$-1$&$-1$&$8$&$8$&$16$&$8$&$16$&$[8]$\\
$139334=2\cdot 233\cdot 13\cdot 23$&$2$&$1$&$1$&$1$&$1$&$-1$&$-1$&$8$&$8$&$16$&$8$&$16$&$[8]$\\
$159310=2\cdot 89\cdot 5\cdot 179$&$2$&$1$&$1$&$1$&$1$&$-1$&$-1$&$8$&$8$&$32$&$8$&$16$&$[16]$\\
\hline
\end{longtable}

\end{example}
\begin{theorem}\label{thm14}  Let $\delta \in \{2p_1, p_2, 2q\}$  such that  $\delta(x\pm1)$   is a square in  $\mathbb{N}$. The order $\#\mathrm{Cl}_2(\mathds{k}^{(1)}_2)=2$  if and only if   the following conditions are satisfied:

\begin{enumerate}[\rm i.]
\item   $ \delta=p_2$ $($resp. $\delta\neq p_2)$ and $(z\pm1)$  $($resp. $\delta(z\pm1))$ is not a square in  $\mathbb{N}$,
\item  at least one of the elements  $\left\{\left(\frac{A}{p_1}\right) ,\left(\frac{2q}{p_1}\right)_4\right\}$  equals  $-1$,
\item $ \left(\frac{p_1}{p_2}\right)_4 \neq \left(\frac{p_2}{p_1}\right)_4 $.
\end{enumerate}
\end{theorem}
\begin{example} Keep notations of Example \ref{9}. For the case: $p_2(x\pm1)$ is a square, $z+1$ and $z-1$ are not squares,  $(\alpha=-1$ or $s=-1)$ and $t_1\neq t_2$.
\begin{longtable}{| c | c | c | c | c |c | c | c | c| c | c | c |c|c|c|c|}
\hline
$d=2p_1.p_2.q  $&$ q_1 $&$ q_2 $&$ q_3 $&$ \alpha $&$ s $&$ t_1 $&$ t_2  $&$ n $&$ n_1 $&$ n_2 $&$ n_3 $&$ q_0$&$ c$\\
\hline
\endfirsthead
\hline
$d=2p_1p_2q  $&$ q_1 $&$ q_2 $&$ q_3 $&$ \alpha $&$ s $&$ t_1 $&$ t_2  $&$ n $&$ n_1 $&$ n_2 $&$ n_3 $&$ q_0$&$ c$\\
\hline
\endhead
$45526=2\cdot 17\cdot 13\cdot 103$&$2$&$1$&$2$&$-1$&$-1$&$-1$&$1$&$8$&$8$&$8$&$8$&$16$&$[6]$\\
$53710=2\cdot 41\cdot 5\cdot 131$&$2$&$1$&$2$&$-1$&$1$&$1$&$-1$&$8$&$8$&$8$&$8$&$16$&$[2]$\\
$56134=2\cdot 17\cdot 13\cdot 127$&$2$&$1$&$2$&$-1$&$1$&$-1$&$1$&$16$&$16$&$16$&$16$&$16$&$[6]$\\
$63438=2\cdot 97\cdot 109\cdot 3$&$2$&$1$&$2$&$-1$&$1$&$1$&$-1$&$8$&$8$&$8$&$8$&$16$&$[2]$\\
\hline
\end{longtable}

\end{example}


\EditInfo{October 15, 2019}{February 28, 2022}{Atilla Berczes}

\end{paper}